\documentclass[num-refs]{wiley-article}
\usepackage{fix-cm}

\usepackage{caption}
\usepackage{setspace}
\usepackage{bm}
\usepackage{stmaryrd}
\usepackage{booktabs}
\usepackage[utf8]{inputenc}
\usepackage{etoolbox}
\usepackage{tabularx}
\usepackage{graphicx}
\usepackage{pdfpages}
\usepackage{incgraph,tikz}
\usepackage[parfill]{parskip}
\usepackage[usestackEOL]{stackengine}
\usepackage{comment}
\usepackage{etoolbox}

\papertype{}
\paperfield{ }

\title{Robust Estimation With Latin Hypercube Sampling: A Central Limit Theorem for Z-estimators}

\author[1,2]{Faouzi Hakimi}

\affil[1]{Aix-Marseille University\\  13284 Marseille, France}
\affil[2]{Institut Mathématique de Toulouse\\  31062 Toulouse, France}

\corraddress{Faouzi Hakimi, Marseille, France}
\corremail{faouzi.hakimi@univ-amu.fr}

\runningauthor{F.Hakimi}

\begin{document}


\begin{frontmatter}

\maketitle

\begin{abstract}

Latin Hypercube Sampling (LHS)  is a widely used stratified sampling method in computer experiments. In this work, we extend existing convergence results for the sample mean under LHS to the broader class of $Z$-estimators — estimators defined as the zeros of a sample mean function. We derive the asymptotic variance of these estimators and demonstrate that it is smaller when using LHS compared to traditional independent and identically distributed  sampling. Furthermore, we establish a Central Limit Theorem for $Z$-estimators under LHS, providing a theoretical foundation for its improved efficiency.

\end{abstract}
\end{frontmatter}

\section{Introduction} \label{Intro}

Latin Hypercube Sampling (LHS), introduced in \cite{LHS}, is a compelling alternative to \emph{independent and identically distributed (i.i.d.)} random sampling for exploring the behavior of complex systems (often treated as black-boxes) through computer experiments \cite{Helton, Helton2, viana2016}.  A significant number of subsequent sampling methods are founded upon LHS (see \cite{DEUTSCH} ,\cite{SHIELDS}, \cite{packham2008latin} or \cite{MONDAL} for instance), underscoring its importance in the field. Consequently, the comprehension and theoretical development of LHS is fundamental in contemporary research. 

\begin{figure}[ht!]
\centering
\includegraphics[width=1.000\textwidth]{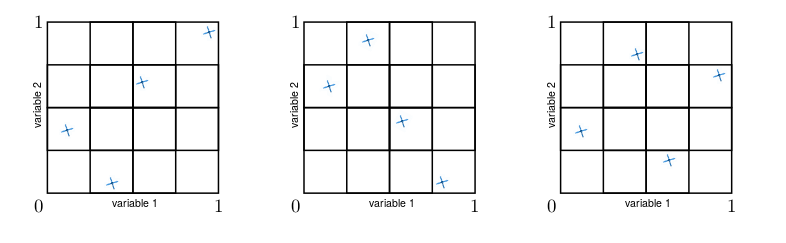}
\caption{Three schematic examples of LHS designs with dimension $d = 2$ and size $n = 4$. \label{LHS2D}}
\end{figure}

To generate an LHS sample of size $n$, the range of each variable is divided into $n$ equally probable intervals. In the case of two variables, the $n$ sample points are then positioned such that there is exactly one sample in each row and each column. Figure \ref{LHS2D} illustrates schematic examples of LHS designs with dimension $d = 2$ and size $n = 4$. The process generalizes naturally to higher dimensions.

Several theoretical results have been established regarding the convergence of estimators under LHS. Most of these results focus on the empirical mean of a measurable function with a finite second-order moment. For instance, it was shown in \cite{Stein} that the asymptotic variance of the sample mean is smaller under LHS compared to classic \emph{i.i.d.} random sampling for such statistics. Additionally, a Central Limit Theorem (CLT) for the empirical mean of bounded functions was proven in \cite{Owen} and later extended to functions with finite third-order moments in \cite{Loh}. 

Subsequent research has extended these results on Latin Hypercube Sampling (LHS) to various methods derived from LHS, such as those accounting for dependencies (e.g., \cite{Aistleitner2013}) and other derivative approaches (e.g, \cite{SHIELDS}). However, the majority of these findings remain focused on the empirical mean.

This paper aims to broaden the scope of the convergence results for empirical mean estimators under classic LHS to encompass the more general class of $Z$-estimators. This class, closely related to the well-established class of $M$-estimators,  includes all estimators that can be expressed as the zeros of an empirical mean function. This work therefore provides novel insights into convergence properties using Latin Hypercube Sampling , which is particularly valuable in computational experiments where variance reduction and efficiency are important to consider.

Most results discussed in this paper are also presented in Chapter 2 of the thesis manuscript \cite{Hakimi}. The paper is organized as follows: Section \ref{LHS} provides a formal definition of Latin Hypercube Sampling along with its key convergence properties. Section \ref{Zestimators} introduces the definitions and relevant properties of $Z$-estimators. Original results concerning the asymptotic normality of $Z$-estimators under LHS are presented in Section \ref{ZLHS}. Finally, an application example is discussed in Section \ref{Ex}.

\newpage
\section{General notations, definitions, and properties of LHS} \label{LHS}

Let us introduce some useful notations regarding this work. We first denote by $\bm{X} = (X_1, \ldots, X_d)$ the vector of $d$ (with $d \in \mathbb{N}^*=\mathbb{N}  \setminus \{0\}$) independent random variables evolving in $\mathcal{\bm{X}} \subset \mathbb{R}^d$. For simplicity and without loss of generality, we assume that the $d$ inputs vary uniformly in $[0,1]$ so we have that, for $j$ in $\llbracket 1,d \rrbracket = \{1,2, \ldots, d$ \} and  $X_j \sim U_{[0,1]}$ . Indeed, one can always work under uniformity and then use the inverse transformation method \cite{Tinverse} to place the support back on the original scale and retrieve the original distribution, as long as the sampling distribution of interest is a product measure (see for instance \cite{Owen} p543 for details). 

Consider a size $n$ ($n \in \mathbb{N}^{*}$) sample of $\bm{X}$ generated using either \emph{i.i.d.} random sampling or Latin Hypercube Sampling (LHS). For either of these methods (denoted here as $METHOD$),  the sample is represented as follows:

\begin{itemize}[wide]
\item The matrix of samples is denoted by $\bm{\mathrm{X}}^{METHOD} = \left( \bm{x}^{(1)},\ldots, \bm{x}^{(n)} \right)^T \in M_{n,d}([0,1])$. Here, $M_{n,d}([0,1])$ denotes the space of matrices of size $n \times d$ with coefficients in $[0,1]$.

\item The $j$th column of $\bm{\mathrm{X}}^{METHOD}$, denoted by $\bm{x}_{j}^{METHOD} = \left(x_{j}^{(1)},\ldots,x_{j}^{(n)}\right)^T$ with $j \in \llbracket 1,d \rrbracket$, represents the effective generated sample of the input $X_j$.
\end{itemize}

Specifically, a sample generated by classic \emph{i.i.d.} random sampling will be denoted by $\bm{\mathrm{X}}^{IID}$, and a sample generated by Latin Hypercube sampling will be denoted by $\bm{\mathrm{X}}^{LHS}$. Similar notations will be used for any quantities estimated with either of these two sampling methods. If no sampling method is mentioned, it means that the results presented do not depend on the sampling method.

We also define the measurable function $g:\mathcal{X} \rightarrow \mathbb{R}^q$ with $q \in \mathbb{N}^{*}$ . This function represents in practice the studied simulation code. We denote by $g(\bm{\mathrm{X}}^{METHOD})= (g (\bm{x}^{(1)}),\ldots,g (\bm{x}^{(n)}) )^T$ the matrix of output samples corresponding to  $\bm{\mathrm{X}}^{METHOD}$.

For $\bm{a} = (a_1, \ldots, a_q) \in \mathbb{R}^q$ with $q \in \mathbb{N}^{*}$, we denote by $||\bm{a}||$ the Euclidean norm of $\bm{a}$ such that $||\bm{a}||^2 = {\sum_{i=1}^q a_i^2}$. Similarly, for any matrix $\bm{A}$ in $M_{q,q}(\mathbb{R})$, we denote by $||\bm{A}||$ the \emph{pseudo} Euclidean norm (Frobenius norm) such that $||\bm{A}||^2 = \sum_{1 \leq i,j \leq q}A_{i,j}^2$. Here, $A_{i,j}$ with ${i,j} \in \llbracket 1,q \rrbracket$ are the components of the matrix $\bm{A}$. A \emph{pseudo} Euclidean norm $||.||$ is finally associated with the tensor space ${T}_{q,q,q}(\mathbb{R})$. This norm is defined, for all $\bm{T}$ in ${T}_{q,q,q}(\mathbb{R})$, by $||\bm{T}||^2 = \sum_{1 \leq i,j,k\leq q}T_{i,j,k}^2$. Here, $T_{i,j,k}$ with ${i,j,k} \in \llbracket 1,q \rrbracket$ are the components of the tensor $\bm{T}$.

We denote by $o(1)$ ("small oh-one") a deterministic sequence that converges to $0$ and $O(1)$  ("big oh-one") a deterministic sequence that is bounded. We denote by $o_p(1)$ ("small oh-P-one") a sequence of random variables that converges in probability to $0$. The expression $O_p(1)$ ("big oh-P-one") denotes a sequence of random variables that is bounded in probability. We recall that a sequence of random variables ${(\bm{W}_n)}_{n \in \mathbb{N}}$ is bounded in probability if, for any scalar $\epsilon >0$, there exist $M$ and $N$ such that, for all $n>N$, $\mathbb{P}(||\bm{W}_n||>M)<\epsilon$ (note that this definition holds in the general case where the norm $||.||$ is not Euclidean).

Finally, a multivariate normal distribution of dimension $q$ ($q \in \mathbb{N}^{*}$) with a mean equal to $\bm{\mu} \in \mathbb{R}^q$ and a covariance matrix equal to $\bm{\Sigma}$ in $M_{q,q}(\mathbb{R})$ is denoted $\mathcal{N}_q(\bm{\mu},\bm{\Sigma})$.


As previously stated, Latin Hypercube Sampling is a statistical method used to generate a near-random sample of parameter values from a multidimensional distribution. To define it formally, we denote, for $d,n \in \mathbb{N}^{*}$:
\begin{enumerate}[wide]
    \item $\bm{\pi}_j=(\pi_{j} (1) \ldots \pi_{j} (n))^T$, $j \in \llbracket 1, d \rrbracket$ as a random permutation of $\llbracket 1, n \rrbracket$, according to the uniform distribution on the set of all possible permutations of $\llbracket 1,n\rrbracket$. Random permutations $(\bm{\pi}_j)_{j \in \llbracket 1, d \rrbracket}$ are assumed to be independent.
    
    \item $\bm{u}_j =(u_{j}^{(1)},\ldots,u_{j}^{(n)})^T, j \in \llbracket 1, d \rrbracket$ as an \emph{i.i.d.} sample of the uniform distribution $U_{[0,1]}$. The samples
$(\bm{u}_j)_{j \in \llbracket 1, d \rrbracket}$ are assumed to be independent.

\end{enumerate}

Random permutations $(\pi_j)_{j \in \llbracket 1, d \rrbracket}$ and samples $(\bm{u}_j)_{j \in \llbracket 1, d \rrbracket}$ are also assumed to be independent. The $n$-sized sampling $x_j^{LHS}$ of the input $X_j, j \in \llbracket 1 ,d \rrbracket$, is then defined as follows:

\begin{equation}
\label{eq:Sample_LHS}
\bm{x}_{j}^{LHS} = \big(x_j^{(1)},\ldots, x_j^{(n)}\big)^T = \left(\frac{1}{n}(\pi_{j} (1) -u_{j}^{(1)}),\ldots, \frac{1}{n}(\pi_{j} (n) -u_{j}^{(n)})\right)^T.
\end{equation}

The corresponding LHS design of dimension $d$ and size $n$ is then $\text{\bf{X}}^{LHS} = (\bm{x}_{1}^{LHS}, \ldots, \bm{x}_{d}^{LHS})$. 


As a result of its stratified nature, the realizations of the LHS design are not \emph{i.i.d.}. However, several results have been indeed established for the convergence of estimators under LHS. Most of them concern the sample mean (first order U-statistics) of measurable functions. For instance, it has been shown in \cite{LHS} that, for any measurable function, this estimator is unbiased:


\begin{proposition}
\label{biais_mean}
Let $g: [0,1]^d \rightarrow \mathbb{R}^q$ with $d,q \in \mathbb{N}^{*}$ be a measurable function such that $\mathbb{E}\big(|| g (\mathbf{X}) || \big)< +\infty$. Denote
$$G^{LHS}_n = \frac{1}{n} \sum_{i=1}^{n} g (\bm{x}^{(i)}),$$
where $\bm{x}^{(i)}, i \in \llbracket 1,n \rrbracket$ is such that $\bm{\mathrm{X}}^{LHS} = \left( \bm{x}^{(1)},\ldots, \bm{x}^{(n)} \right)^T$ with $\bm{\mathrm{X}}^{LHS}$ being defined using Equation \eqref{eq:Sample_LHS}. Then, $G^{LHS}_n$ is an unbiased estimator of $G = \mathbb{E} \big( g(\bm{X}) \big)$.
\end{proposition}
Similarly to Proposition \ref{biais_mean}, we denote by $G^{IID}_n$ the classic sample mean of an \emph{i.i.d.} design: $G^{IID}_n = \frac{1}{n} \sum_{i=1}^{n} g (\bm{x}^{(i)}),$ with $\bm{x}^{(i)}, i \in \llbracket 1,n \rrbracket$  being a sample of an \emph{i.i.d.} design $\bm{\mathrm{X}}^{IID} \in M_{n,d}([0,1])$. 

A second interesting characteristic of mean value estimators under LHS is their variance. Indeed, Stein \cite{Stein} showed that if $g$ is a real-valued  function such that $\mathbb{E} \big(g^2 (\mathbf{X}) \big) <+\infty$, then $\mathbb{V}\text{ar}(G^{LHS}_n)$ is always asymptotically smaller than $\mathbb{V}\text{ar}(G^{IID}_n)$. This result is generalized  to multidimensional functions by Loh in \cite{Loh}. Proposition  \ref{var_mean} summarizes the main results regarding the covariance matrix of $G^{LHS}_n$:

\begin{proposition}
\label{var_mean}
Let $g: [0,1]^d \rightarrow \mathbb{R}^q$  ($d,q \in \mathbb{{N}}^{*}$) be a measurable function with $\mathbb{E}\big( || g (\mathbf{X}) ||^2 \big) < +\infty$. Let  $ \bm{\Sigma}_{G^{IID}_n}, \bm{\Sigma}_{G^{LHS}_n} \in M_{q,q}(\mathbb{R})$
be the covariance matrices of $G^{IID}_n$ and $G^{LHS}_n$ respectively, with $\bm{\Sigma}_{G^{IID}_n}= \frac{1}{n}\mathbb{E}\bigg(\big(g(\bm{X})-G\big)\big(g(\bm{X})-G\big)^T \bigg)$.

We also define, for $\bm{x}=(x_1,\ldots, x_d) \in [0,1]^d$:

\begin{itemize}[wide]

\item $g_{-j}(x_j) = \int_{[0,1]^{d-1}}[g(\bm{x})-G] \underset{1 \leq k \leq  d, k \neq j}{\prod}{dx_k} = \mathbb{E} \big( g(\bm{X})- G|X_j \big) $ with $j \in \llbracket 1,d \rrbracket$.

\item $g_{rem}(\bm{x}) = g(\bm{x}) -  G - \sum_{j=1}^d g_{-j}(x_j).$ 

\item $\bm{R}_g=\int_{[0,1]^d}g_{rem}(\bm{x})g_{rem}(\bm{x})^T\bm{dx}$. 

\end{itemize}

Then we have:

\begin{itemize}[wide]

\item  $\bm{\Sigma}_{G^{LHS}_n} =  \frac{1}{n}\bm{R}_g+\frac{1}{n} o(1)$.

\item $\bm{\Sigma}_{G^{IID}_n} =  \frac{1}{n}\bm{R}_g + \frac{1}{n}\sum_{j=1}^d \int_{[0,1]} g_{-j}(x_j)g_{-j}(x_j)^T dx_j$.

\end{itemize}

We therefore have that $\bm{\Sigma}_{G^{IID}_n} -\bm{\Sigma}_{G^{LHS}_n}$ is asymptotically positive semi-definite, that is,

$\forall \xi \in \mathbb{R}^d, \quad \underset{n\rightarrow +\infty} {\lim}n \xi^T (\bm{\Sigma}_{G^{IID}_n} -\bm{\Sigma}_{G^{LHS}_n})\xi\geq
\sum_{j=1}^d \int_{[0,1]} \xi^T g_{-j}(x_j)g_{-j}(x_j)^T \xi dx_j
\geq0.$
\end{proposition}

Since  $G^{IID}_n$ converges in quadratic mean to $G$ and that  $G^{LHS}_n$ is an unbiased estimator of $G$ (as established in Proposition \ref{biais_mean}), we can conclude that  $G^{LHS}_n$ also converges in quadratic mean to $G$: $\lim\limits_{n\rightarrow +\infty} \mathbb{E} \big( ||G^{LHS}_n- G|| ^2 \big) = 0$. Consequently, $G^{LHS}_n$  converges in probability to $G$.

In addition, Owen \cite{Owen} showed a Central Limit Theorem (CLT) for this class of estimators under LHS when the model function $g$ is bounded. This was generalized to any function with finite third moment in \cite{Loh}:

\begin{theorem}
\label{TCL_mean}
In the framework of Proposition \ref{var_mean}, let $g: [0,1]^d \rightarrow \mathbb{R}^q$ ($d,q \in \mathbb{N}^{*}$) be a measurable function with $\mathbb{E} \big( ||g(\mathbf{X})||^3 \big) < +\infty$. Then, assuming that $\bm{R}_g$ is non-singular, we have that $\sqrt{n}(G_n^{LHS}-G)$ tends in distribution to $\mathcal{N}_q(0,\bm{R}_g)$ as ${n \to +\infty}$. 

\end{theorem}

\section{ Definitions and properties on \textit{Z}-estimators  } \label{Zestimators}

The primary objective of this work is to extend the convergence results under Latin Hypercube Sampling (LHS) to the class of $Z$-estimators. The $Z$-estimator class is intimately related to the well-established class of $M$-estimators, yet it offers a distinct formulation. Specifically, $Z$-estimators are defined as solutions to a set of estimating equations, which can be viewed as a generalization of the optimization problem associated with $M$-estimators. This formulation provides a flexible framework for parameter estimation, encompassing a wide range of statistical models and inference procedures. For a comprehensive discussion on these topics, we refer the reader to \cite{Robust} and \cite{VanderVaart}. Significant theoretical work has been conducted on this class of estimators, which remains a major research topic today (see, for example, \cite{Bouzebda2022} or \cite{Zhan2002}).

More formally, let $\bm{\text{X}}= (\bm{x}^{(1)},\ldots,\bm{x}^{(n)} )^T$ be the vector of $n$ realizations of a random vector $\bm{X}$ evolving in  $\bm{\mathcal{X}} \subset \mathbb{R}^d$, with $n,d \in \mathbb{N}^{*}$. Its law is parameterized by a vector $\bm{\theta} \in \Theta \subset \mathbb{R}^q, q \in \mathbb{N}^{*}$. 

For $\bm{x} \in \mathcal{X}, \theta \in \Theta$, let $(\bm{x}, \bm{\theta}) \rightarrow \psi_{\bm{\theta}}(\bm{x}) \in \mathbb{R}^q$  be a known measurable function such that $\psi_{\bm{\theta}}(\bm{x}) = ( \psi_{\theta_1} (\bm{x}), \ldots, \psi_{\theta_q}(\bm{x}) )^T$. We also define the empirical mean of this function $(\bm{\text{X}}, \bm{\theta}) \rightarrow \Psi_n(\bm{\theta})\in \mathbb{R}^q$ such that $\Psi_n(\bm{\theta}) = \frac{1}{n}\sum_{i=1}^n \psi_{\bm{\theta}}(\bm{x}^{(i)})$.  

The  $Z$-estimator $\hat{\bm{\theta}}_n = \hat{\bm{\theta}}_n(\bm{x}^{(1)},\ldots,\bm{x}^{(n)} ) \in \Theta$ associated with $\psi_{\bm{\theta}}(\bm{x})$ corresponds to the solution of the following vectorial equation: 
\begin{equation}
    \Psi_n(\bm{\theta}) =0.
\label{eq:Zestimator}
\end{equation}

Many known estimators can be defined as $Z$-estimators. For instance, let $\bm{X}$ have a distribution function $f_{\bm{\theta}}$ with a continuous first derivative in ${\bm{\theta}} \in \Theta$ and a separable log-likelihood function. In this case, the maximum likelihood estimator of $\bm{\theta}$ can be written as a $Z$-estimator as defined by \ref{eq:Zestimator} with, for $\bm{x} \in \mathbb{R}^d, d \in \mathbb{N}^{*}$, $\psi_{\bm{\theta}}(\bm{x}) = (\frac{\partial \log(f_{\bm{\theta}}(\bm{x}))}{\partial \theta_1},\ldots, \frac{\partial \log(f_{\bm{\theta}}(\bm{x}))}{\partial \theta_q})^T$. The $Z$-estimator framework also encompasses Generalized Method of Moments (GMM) estimators  \cite{GMM}, where $\psi_\theta(x)$ comprises moment conditions $g(x, \theta)$ equating sample and theoretical moments. Quantile regression also fits within this framework, with $\psi_\theta(x, y) = [I(y \leq x^T \theta) - \tau] x$, where $\tau$ is the quantile of interest, and $I(\cdot)$ is the indicator function. This estimates conditional quantiles rather than the mean. These examples highlight the versatility of the $Z$-estimator framework, where the key is specifying an appropriate $\psi$ function reflecting the estimation objectives.

The first useful properties regarding $Z$-estimators concern the link between the consistency of $\Psi_n(\bm{\theta})$ and the consistency of $\hat{\bm{\theta}}_n$. For instance, in \cite{Duflo}, one can find assumptions for which the consistency of $\bm{\hat{\theta}}_n$ is ensured:

\begin{proposition}
\label{convergence_2}

Let $\Theta$ be a compact subset of $\mathbb{R}^q$ with $q \in \mathbb{N}^{*}$. Let also assume that the following hypotheses are true, for any $\bm{\theta} \in \Theta$ and $n \in \mathbb{N}^{*}$:

\begin{itemize}[wide]

\item the functions $\bm{\theta} \rightarrow \Psi_n(\bm{\theta}) $ and $\bm{\theta} \rightarrow \Psi(\bm{\theta})$ are continuous measurable functions of $\bm{\theta} \in \Theta$ evolving in $\mathbb{R}^q$;

\item each function $\Psi_n(\bm{\theta})$ has exactly one zero $\hat{\bm{\theta}}_n \in \Theta$;

\item $\Psi_n(\bm{\theta})$ converges to  $\Psi(\bm{\theta})$ in probability;

\item $\Psi(\bm{\theta})$ vanishes only at $\bm{\theta}_0$ with $\bm{\theta_0} \in \Theta$;

\item denoting, for $\eta \geq 0$, $w_n(\eta)=\sup \{|| \Psi_n(\bm{\theta}_1) - \Psi_n(\bm{\theta}_2) ||; ||\bm{\theta}_1-\bm{\theta}_2 ||\leq \eta, \bm{\theta_1}, \bm{\theta_2} \in \Theta \}$;
there exists two sequences $(\eta_k)$ and $(\epsilon_k)$ both decreasing to $0$ such that, for all $k \in \mathbb{N},$ $\mathbb{P}(w_n(\eta_k)>\epsilon_k) \xrightarrow[n \rightarrow +\infty]{}0$.

\end{itemize}
Then $\bm{\hat{\theta}}_n$ is a consistent estimator of $\bm{\theta}_0$, that is  $\bm{\hat{\theta}}_n \xrightarrow[n \rightarrow + \infty]{p} \bm{\theta}_0$.

\end{proposition}

The assumption on $w_n(\eta)$ seems difficult to grasp at first glance. However, as mentioned in \cite{Duflo}, if we find a function $\phi$ from $\mathbb{R}_+$ to $\mathbb{R}$ such that  $\underset{\eta \rightarrow 0^+}{\lim}\phi(\eta)=0$, this assumption on $w_n$ can be obtained through: $\mathbb{P}(w_n(\eta)\geq 2 \phi(\eta)) \xrightarrow[n\rightarrow +\infty]{}0$ for each $\eta\geq0$. For instance, $w_n(\eta) \xrightarrow[n\rightarrow +\infty]{}\phi(\eta)$, or $\underset{n\rightarrow +\infty}{\lim} w_n(\eta)\leq \phi(\eta)$ give both sufficient conditions. Note that Proposition \ref{convergence_2} is general and does not mention any sampling scheme. 

In addition to these convergence properties, several Central Limit Theorems for $Z$-estimators have been proved. Here we give one of them, proposed in \cite{VanderVaart}. Theorem \ref{TCL} relies on the so-called \emph{classic conditions}, formulated  to mathematically tighten the informal derivation of the asymptotic normality of maximum likelihood proposed by  \cite{Fisher_likelihood}. These conditions are stringent, but they are simple. They lead to a simple proof of the Central Limit Theorem. This simplicity will allow us to adapt this theorem to the LHS case.

In particular, a needed assumption for the application of this theorem concerns the existence of a first and  a second order derivatives in $\bm{\theta}$ for  $\psi_{\bm{\theta}}$. Let us introduce these terms.

For any $\bm{\theta} \in \Theta$ and for any $\bm{x} \in \mathcal{X}$, let $(\bm{x}, \bm{\theta}) \rightarrow \dot{\psi}_{\bm{\theta}}(\bm{x})$ be the first order partial derivative of $\psi_{\bm{\theta}} \in \mathbb{R}^q$, assuming it exists. This first order partial derivative is evolving in $M_{q,q}(\mathbb{R})$ . Its components are  $\dot{\psi}_{\bm{\theta}_{j,k}}= \frac{\partial \psi_{\bm{\theta}_j}}{\partial \theta_k }$ with $j,k \in \llbracket1,q \rrbracket$.

Similarly, for any $\bm{\theta} \in \Theta$ and for any $\bm{x} \in \mathcal{X}$, let $(\bm{x}, \bm{\theta}) \rightarrow \ddot{\psi}_{\bm{\theta}}(\bm{x})$ be the second order partial derivative of $\psi_{\bm{\theta}} \in \mathbb{R}^q$, assuming it exists. This second order partial derivative is evolving in $T_{q,q,q}(\mathbb{R})$ . Its components are  $\ddot{\psi}_{\bm{\theta}_{j,k,l}}= \frac{\partial^2 \psi_{\bm{\theta}_j}}{\partial \theta_k \partial \theta_l}$, with $j,k,l \in \llbracket1,q \rrbracket$.


\begin{theorem}
\label{TCL}

Let $\Theta$ be an open subset of an Euclidean space of dimension $q, q \in {\mathbb{N}^{*}}$ and let $\mathcal{X}$ be a subspace of $\mathbb{R}^d, d \in {\mathbb{N}^{*}}$.  Assume that, for all $\bm{\theta}$ in $\Theta$ and for all $\bm{x}$ in $\mathcal{X}$, the function $(\bm{x},\bm{\theta}) \rightarrow \psi_{\bm{\theta}}(\bm{x})$ evolving in $\mathbb{R}^q$ is twice continuously differentiable in $\bm{\theta}$.

Let $(\bm{x}, \bm{\theta}) \rightarrow \dot{\psi}_{\bm{\theta}}(\bm{x}) \in M_{q,q}(\mathbb{R})$ and $(\bm{x}, \bm{\theta}) \rightarrow \ddot{\psi}_{\bm{\theta}}(\bm{x}) \in T_{q,q,q}(\mathbb{R})$ denote the first and second-order derivatives of $\psi_{\bm{\theta}}$, respectively.

Let $\bm{\text{X}}^{IID} = (\bm{x}^{(1)},\ldots,\bm{x}^{(n)} )^T$ be the vector of i.i.d. realizations of a random variable $\bm{X}= (X_1,\ldots X_d)$ evolving in $\mathcal{X}$.

Suppose also that the following assumptions are fulfilled:

\begin{enumerate}[wide]
\item $\Psi_n^{IID}(\hat{\bm{\theta}}_n^{IID}) = \frac{1}{n}\sum_{i=1}^n \psi_{\hat{\bm{\theta}}_n^{IID}}(\bm{x}^{(i)}) =  0, \forall n \in \mathbb{N}^*$;
\item there exists a unique $\bm{\theta}_0$ in $\Theta$ such that $\mathbb{E}(\psi_{\bm{\theta}_0}(\bm{X}))= \Psi(\bm{\theta}_0) = 0$ with $\bm{\theta}_0$ in $\Theta$;
\item $\mathbb{E}(||\psi_{\bm{\theta}_0}(\bm{X})||^2)<+\infty$;
\item $\mathbb{E}(\dot{\psi}_{\theta_0}(\bm{X}))$ exists and is non-singular;
\item For any $ \bm{x} \in \mathcal{X}$ and  for any $\bm{\theta}$ in the neighborhood of $\bm{\theta}_0$, the function $(\bm{x},\bm{\theta}) \rightarrow \ddot{\psi}_{\bm{\theta}}(\bm{x}) \in T_{q,q,q}(\mathbb{R})$ is dominated,  in norm, by a fixed integrable function $\bm{x} \rightarrow \ddot{\psi}(\bm{x}) \in T_{q,q,q}(\mathbb{R})$.
\end{enumerate}

Then, if  $\hat{\bm{\theta}}_n^{IID}$ is a consistent  estimator of $\bm{\theta}_0$ , we have:

\begin{equation}
    (\hat{\bm{\theta}}^{IID}_n-\bm{\theta}_0) = -[\mathbb{E}(\dot{\psi}_{\bm{\theta}_0}(\bm{X}))]^{-1}\frac{1}{{n}}\sum_{i=1}^n \psi_{\bm{\theta}_0}(\bm{x}^{(i)}) + \frac{1}{\sqrt{n}}o_p(1).
\end{equation}

Moreover, we have  that the sequence $\sqrt{n}(\hat{\bm{\theta}}^{IID}_n-\bm{\theta}_0)$ tends in distribution to  \\ $\mathcal{N}_q(0,[\mathbb{E}(\dot{\psi}_{\bm{\theta}_0}(\bm{X}))]^{-1}
\mathbb{E}\big(\psi_{\bm{\theta_0}}(\bm{X})\psi_{\bm{\theta}_0}(\bm{X})^T\big)
[\mathbb{E}(\dot{\psi}_{\bm{\theta}_0}(\bm{X}))]^{-T})$ as $n \to + \infty$.

\end{theorem}

For the following, it is important to note that we have  $\bm{\Sigma}_{\Psi_n^{IID}({\bm{\theta}}_0)} = \frac{1}{n}\mathbb{E}\big(\psi_{\bm{\theta_0}}(\bm{X})\psi_{\bm{\theta}_0}(\bm{X})^T\big) $ with $ \bm{\Sigma}_{\Psi_n^{IID}({\bm{\theta}}_0)} \in M_{q,q}(\mathbb{R})$ being the covariance matrix of $\Psi_n^{IID}({\bm{\theta}}_0)$. It is also important to remark that among the results presented in this section, only Theorem \ref{TCL} requires the specific use of an \emph{i.i.d.} sample, since its proof relies on the classical Central Limit Theorem (CLT) \cite{ClassicTCL}.

While Theorem \ref{TCL} assumes $\Theta$ is open, one can modify this assumption to consider $\Theta$  as the interior of a compact set. This allows us to maintain the differentiability conditions required for asymptotic normality while preserving the compactness needed for Proposition \ref{convergence_2} concerning consistency.

\section{{\textit{Z}}-estimators under LHS} \label{ZLHS}

In this section, we extend the convergence properties of $Z$-estimators to LHS designs. The idea is to combine all the above properties. Indeed, one can first notice that the $Z$-function $\Psi_n(\bm{\theta})$ is the empirical mean of $\psi_{\bm{\theta}}$. Now, as mentioned in Section \ref{LHS}, the convergence of this type of statistic under LHS holds. We use that here to show a Central Limit Theorem for $Z$-estimators under LHS.

As in Section \ref{LHS},  let, for any $\bm{\theta} \in \Theta$ and $\bm{\text{X}}^{LHS}, \bm{\text{X}}^{IID} \in M_{n,d}([0,1])$, $\bm{\Sigma}_{\Psi_n^{IID}({\bm{\theta}})} ,\bm{\Sigma}_{\Psi_n^{LHS}({\bm{\theta}})} \in M_{q,q}(\mathbb{R})$ be the covariance matrices of  $\Psi_n^{IID}(\bm{\theta})$ and $\Psi_n^{LHS}(\bm{\theta})$ respectively. Let us now give some noteworthy convergence properties on $\Psi_n^{LHS}(\bm{\theta})$.

\begin{proposition}
\label{PsiConv_1}

Let $\Theta$ be a compact subset of  $\mathbb{R}^q$ and $\mathcal{X}= [0,1]^d$ ($q,d$ in ${\mathbb{N}^{*}}$). Let $\bm{\text{X}}^{LHS} = (\bm{x}^{(1)},\ldots,\bm{x}^{(n)} )^T$ be the vector of LHS realizations of a random variable $\bm{X}= (X_1,\ldots X_d)$ evolving in $\mathcal{X}$ such that $\bm{X} \sim U_{[0,1]^d}$. Assume also that, for all $\bm{\theta} \in \Theta$ and $\bm{x} \in \mathcal{X}$, the function $(\bm{x}, \bm{\theta}) \rightarrow \psi_{\bm{\theta}}(\bm{x})$ is measurable regarding $\bm{x}$. We then have the following properties on $\Psi_n^{LHS}(\bm{\theta}) = \frac{1}{n}\sum_{i=1}^n \psi_{\bm{\theta}}(\bm{x}^{(i)})$:

\begin{enumerate}[wide]
\item  If, for all $\bm{\theta} \in \Theta$,  $\mathbb{E}(|| \psi_{\bm{\theta}}(\bm{X})||) <+\infty$, $\Psi_n^{LHS}(\bm{\theta})$ is an unbiased estimator of $\Psi(\bm{\theta}) =  \mathbb{E}(\psi_{\bm{\theta}}(\bm{X}))$.

\item  If, for all $\bm{\theta} \in \Theta$,  $\mathbb{E}(|| \psi_{\bm{\theta}}(\bm{X})||^2) <+\infty$, we also have:

$\bm{\Sigma}_{\Psi_n^{LHS}({\bm{\theta}})}=  \frac{1}{n}\int_{[0,1]^d}{\psi_{\bm{\theta}_{rem}}}(\bm{x})\psi_{\bm{\theta}_{rem}}(\bm{x})^T\bm{dx} + \frac{1}{n}o(1)$, with $\psi_{\bm{\theta}_{rem}}$ being defined as in Proposition \ref{var_mean}. 

Moreover, we have that $\bm{\Sigma}_{\Psi_n^{IID}({\bm{\theta}})}- \bm{\Sigma}_{\Psi_n^{LHS}({\bm{\theta}})}$ is asymptotically positive semi-definite and that $\Psi_n^{LHS}(\bm{\theta})$ converges in quadratic mean to $\Psi(\bm{\theta})$. In other words, we have  $\lim\limits_{n\rightarrow +\infty} \mathbb{E} \big(||\Psi_n^{LHS}(\bm{\theta})- \Psi(\bm{\theta})||^2 \big) = 0$.

\item If, for all $\bm{\theta} \in \Theta$,  $\mathbb{E}(|| \psi_{\bm{\theta}}(\bm{X})||^3) <+\infty$ and if $\bm{R}_{\psi_{\bm{\theta}}} = \int_{[0,1]^d}{\psi_{\bm{\theta}_{rem}}}(\bm{x})\psi_{\bm{\theta}_{rem}}(\bm{x})^T\bm{dx}$ is non-singular, we have that  $\sqrt{n}(\Psi_n^{LHS}-\Psi(\bm{\theta}))$ tends in distribution to $\mathcal{N}_q(0,\bm{R}_{\psi_{\bm{\theta}}})$ as $n \to +\infty$.
\end{enumerate}

\end{proposition}

\emph{Proof.} Let us show these properties one by one:
\begin{enumerate}[wide]
    \item Since, for all $\bm{\theta} \in \Theta$ and $\bm{x} \in \mathcal{X}$, the function $(\bm{x}, \bm{\theta}) \rightarrow \psi_{\bm{\theta}}(\bm{x})$ is measurable regarding $\bm{x} \in \mathcal{X}$ and $\mathbb{E}(|| \psi_{\bm{\theta}}(\bm{X})||) <+\infty$, $\Psi_n^{LHS}(\bm{\theta})$ is an unbiased estimator of $\Psi(\bm{\theta})$ by Proposition \ref{biais_mean}.
    \item This is a direct consequence of Proposition \ref{var_mean}.
    \item This is a direct consequence of Theorem \ref{TCL_mean}. $\blacksquare$
\end{enumerate}
 
All these properties on $\Psi_n^{LHS}(\bm{\theta})$ allow to show that $\hat{\bm{\theta}}_n^{LHS}$ is a consistent estimator of $\bm{\theta}_0$. Indeed, the assertion 2 of Proposition \ref{PsiConv_1} ensures the convergence in probability of $\Psi_n^{LHS}(\bm{\theta})$ to $\Psi(\bm{\theta})$. As mentioned before,  Proposition \ref{convergence_2} does not impose any other conditions on the sampling scheme. We therefore have, under the conditions of applicationof this proposition, that $\hat{\bm{\theta}}_n^{LHS}$ converges to $\bm{\theta}_0$ in probability. Let us now establish a Central Limit Theorem for $Z$-estimators under LHS.

\begin{theorem}
\label{ZLHS_TCL}

Let $\Theta$ be the interior of a compact subset of $\mathbb{R}^q$, $q \in \mathbb{N}^*$, and $\mathcal{X} = [0,1]^d$, $d \in \mathbb{N}^*$. For all $\bm{\theta} \in \Theta$ and $\bm{x} \in \mathcal{X}$, assume $(\bm{x}, \bm{\theta}) \rightarrow \psi_{\bm{\theta}}(\bm{x})$, where $\psi_{\bm{\theta}} = (\psi_{\theta_1}(\bm{x}), \ldots, \psi_{\theta_q}(\bm{x}))^T \in \mathbb{R}^q$, is twice continuously differentiable in $\bm{\theta}$. 

Let $(\bm{x}, \bm{\theta}) \rightarrow \dot{\psi}_{\bm{\theta}}(\bm{x}) \in M_{q,q}(\mathbb{R})$ and $(\bm{x}, \bm{\theta}) \rightarrow \ddot{\psi}_{\bm{\theta}}(\bm{x}) \in T_{q,q,q}(\mathbb{R})$ denote the first and second-order derivatives of $\psi_{\bm{\theta}}$, respectively. 

For any $n \in \mathbb{N}^{*}$, let $\bm{\text{X}}^{LHS} = (\bm{x}^{(1)},\ldots,\bm{x}^{(n)})^T$ be LHS realizations of $\bm{X} \sim U_{[0,1]^d}$ of size $n$. 
Suppose also that the following hypotheses are fulfilled:
\begin{enumerate}
\item For any $n \in \mathbb{N}^{*}, \Psi_n(\bm{\hat{\theta}}_n^{LHS}) = \frac{1}{n}\sum_{i=1}^n \psi_{\bm{\hat{\theta}}_n^{LHS}}(\bm{x}^{(i)}) = 0$ ;
\item There is $\bm{\theta}_0 \in \Theta$ such that $\mathbb{E}(\psi_{\bm{\theta}_0}(\bm{X})) = \Psi(\bm{\theta}_0) = 0$ ;
\item $\mathbb{E}(\|\psi_{\bm{\theta}_0}(\bm{X})\|^2) < + \infty$ ;
\item $\mathbb{E}(\dot{\psi}_{\bm{\theta}_0}(\bm{X}))$ is non-singular and $\mathbb{E}(\|\dot{\psi}_{\bm{\theta}_0}(\bm{X})\|^2) < + \infty$ ;
\item  There is an integrable function $\bm{x} \rightarrow \ddot{\psi}(\bm{x})  \in T_{q,q,q}(\mathbb{R}), \bm{x} \in \mathcal{X}$, such that $\|\ddot{\psi}_{\bm{\theta}}(\bm{x})\| \leq \|\ddot{\psi}(\bm{x})\|$ and $\mathbb{E}(\|\ddot{\psi}(\bm{X})\|^2) < + \infty$ for all $\bm{\theta}$ in the neighborhood of  $\bm{\theta}_0$.
\end{enumerate}
If $\hat{\bm{\theta}}_n^{LHS}$ is a consistent estimator of $\bm{\theta}_0$, then:
\begin{equation}
\label{eqn:ZLHS_TCL}
\hat{\bm{\theta}}_n^{LHS} - \bm{\theta}_0 = -[\mathbb{E}(\dot{\psi}_{\bm{\theta}_0}(\bm{X}))]^{-1}\frac{1}{n}\sum_{i=1}^n \psi_{\bm{\theta}_0}(\bm{x}^{(i)}) + \frac{1}{\sqrt{n}}o_p(1).
\end{equation}
Moreover, let $\bm{\Sigma}_{\bm{\hat{\theta}}_n^{LHS}} \in M_{q,q}(\mathbb{R})$ be the covariance matrix of $\bm{\hat{\theta}}_n^{LHS}$. Then:
\begin{equation}
\bm{\Sigma}_{\bm{\hat{\theta}}_n^{LHS}} = [\mathbb{E}(\dot{\psi}_{\bm{\theta}_0}(\bm{X}))]^{-1}
\bm{\Sigma}_{\Psi_n^{LHS}({\bm{\theta}}_0)}[\mathbb{E}(\dot{\psi}_{\bm{\theta}_0}(\bm{X}))]^{-T} + \frac{1}{n}o(1),
\end{equation}
where
\begin{equation}
    \bm{\Sigma}_{\Psi_n^{LHS}({\bm{\theta}}_0)} = \frac{1}{n}\bm{R}_{\psi_{\bm{\theta}_0}} + \frac{1}{n}o(1).\end{equation}

with $\bm{R}_{\psi_{\bm{\theta}_0}} = \int_{[0,1]^d}{\psi_{{\bm{\theta}_0}_{rem}}}(\bm{x})\psi_{{\bm{\theta}_0}_{rem}}(\bm{x})^T\bm{dx}$. Furthermore, $\bm{\Sigma}_{\bm{\hat{\theta}}_n^{IID}} - \bm{\Sigma}_{\bm{\hat{\theta}}_n^{LHS}}$ is asymptotically positive semi-definite ($\bm{\Sigma}_{\bm{\hat{\theta}}_n^{IID}} \in M_{q,q}(\mathbb{R})$ corresponding to the covariance matrice of $\bm{\hat{\theta}}_n^{IID}$). \\
Finally, if $\mathbb{E}(\|\psi_{\bm{\theta_0}}(\bm{X})\|^3) < +\infty$ and $\bm{R}_{\psi_{\bm{\theta}_0}}$ is non-singular, then $\sqrt{n}(\bm{\hat{\theta}}^{LHS} - \bm{\theta}_0)$ tends in distribution to \\ $\mathcal{N}(0, [\mathbb{E}(\dot{\psi}_{\bm{\theta}_0}(\bm{X}))]^{-1}\bm{R}_{\psi_{\bm{\theta}_0}}[\mathbb{E}(\dot{\psi}_{\bm{\theta}_0}(\bm{X}))]^{-T})$ as ${n \rightarrow + \infty}$.

\end{theorem}

\emph{Proof.}  The proof follows the reasoning in \cite{VanderVaart} for Theorem \ref{TCL}.

By Taylor's Theorem,  as $\Psi_n(.)$ is continuous  and twice differentiable in $\bm{\theta}$,  $\exists  \tilde{\bm{\theta}}_n^{LHS}$ between $\bm{\theta}_0$ and $\bm{\hat{\theta}}_n^{LHS}$ such that:
\begin{equation}
\begin{split}
\Psi_n^{LHS}(\hat{\bm{\theta}}_n^{LHS}) = 0  &= \Psi_n^{LHS}(\bm{\theta}_0) + \dot{\Psi}_n^{LHS}(\bm{\theta}_0)(\hat{\bm{\theta}}_n^{LHS} - \bm{\theta}_0) \\
&\quad + \frac{1}{2}(\hat{\bm{\theta}}_n^{LHS} - \bm{\theta}_0)^T \ddot{\Psi}_n^{LHS}(\tilde{\bm{\theta}}_n^{LHS})(\hat{\bm{\theta}}_n^{LHS} - \bm{\theta}_0).
\end{split}
\end{equation}
Since $\mathbb{E}(\|\psi_{\bm{\theta}_0}(\bm{X})\|^2) < +\infty$, Proposition \ref{PsiConv_1} implies:
\begin{equation}
\Psi_n^{LHS}(\bm{\theta}_0) = \frac{1}{n}\sum_{i=1}^n \psi_{\bm{\theta_0}}(\bm{x}^{(i)}) \xrightarrow[n\rightarrow +\infty]{p} \mathbb{E}(\psi_{\bm{\theta}_0}(\bm{X})) = 0.
\end{equation}
Now, let  $\dot{\Psi}_n^{LHS}(\bm{\theta})  = \frac{1}{n}\sum_{i=1}^n \dot{\psi}_{\bm{\theta}}(\bm{x}^{(i)})$ be the empirical mean over $\bm{x}$ of the matrix function  $(\bm{x}, \bm{\theta}) \rightarrow \dot{\psi}_{\bm{\theta}}(\bm{x})$, with $\ddot{\psi}_{\bm{\theta}}(\bm{x}) \in {M}_{q,q}(\mathbb{R})$.  Similarly,  $\dot{\Psi}_n^{LHS}(\bm{\theta})  \xrightarrow[n\rightarrow + \infty]{p} \mathbb{E}(\dot{\psi}_{\bm{\theta}_0}(\bm{X}))$, which is non-singular by assumption.

Let also $\ddot{\Psi}_n^{LHS}(\bm{\theta})  = \frac{1}{n}\sum_{i=1}^n \ddot{\psi}_{\bm{\theta}}(\bm{x}^{(i)})$ be the empirical mean over $\bm{x}$  of the tensor function $(\bm{x}, \bm{\theta}) \rightarrow \ddot{\psi}_{\bm{\theta}}(\bm{x})$, with $\ddot{\psi}_{\bm{\theta}}(\bm{x}) \in {T}_{q,q,q}(\mathbb{R})$.

For $\ddot{\Psi}_n^{LHS}(\bm{\theta})$, let $\mathcal{B}$ be a ball around $\bm{\theta}_0$ where $\|\ddot{\psi}_{\bm{\theta}}\| \leq \|\ddot{\psi}\|$ with $\bm{x} \rightarrow \ddot{\psi}(\bm{x})  \in T_{q,q,q}(\mathbb{R})$ being an integrable function and $\mathbb{E}(\|\ddot{\psi}(\bm{X})\|^2) < + \infty$  (this ball exists by assumption). Since $\hat{\bm{\theta}}_n^{LHS} \xrightarrow[n\rightarrow +\infty]{p} \bm{\theta}_0$, we have $\mathbb{P}(\tilde{\bm{\theta}}_n^{LHS} \in \mathcal{B}) \rightarrow 1$. For $\tilde{\bm{\theta}}_n^{LHS} \in \mathcal{B}$:
\begin{equation}
\|\ddot{\Psi}_n^{LHS}(\tilde{\bm{\theta}}_n^{LHS})\| \leq \frac{1}{n} \sum_{i=1}^n \|\ddot{\psi}(\bm{x}^{(i)})\|.
\end{equation}
The right-hand side converges to a finite value by Proposition \ref{var_mean}, implying the same for the left-hand side.

Rewriting the Taylor expansion:

\begin{equation}
-\Psi_n^{LHS}(\bm{\theta}_0) = \left(\mathbb{E}(\dot{\psi}_{\bm{\theta}_0}(\bm{X})) + o_p(1) + \frac{1}{2}(\hat{\bm{\theta}}_n^{LHS} - \bm{\theta}_0)^T O_p(1)\right)(\hat{\bm{\theta}}_n^{LHS} - \bm{\theta}_0).
\end{equation}
As $\hat{\bm{\theta}}_n^{LHS} \xrightarrow[n\rightarrow +\infty]{p} \bm{\theta}_0$, we have:
\begin{equation}
-\Psi_n^{LHS}(\bm{\theta}_0) = \left(\mathbb{E}(\dot{\psi}_{\bm{\theta}_0}(\bm{X})) + o_p(1)\right)(\hat{\bm{\theta}}_n^{LHS} - \bm{\theta}_0).
\end{equation}
This yields equation \eqref{eqn:ZLHS_TCL}, since $\mathbb{E}(\dot{\psi}_{\bm{\theta}_0})$ is non-singular and $\Psi_n^{LHS}(\bm{\theta}_0) = \frac{1}{\sqrt{n}}O_p(1)$ asymptotically. 

We also have $\bm{\Sigma}_{\bm{\hat{\theta}}_n^{LHS}} = [\mathbb{E}(\dot{\psi}_{\bm{\theta}_0}(\bm{X}))]^{-1}
\bm{\Sigma}_{\Psi_n^{LHS}({\bm{\theta}}_0)}[\mathbb{E}(\dot{\psi}_{\bm{\theta}_0}(\bm{X}))]^{-T} + \frac{1}{n}o(1).$ 

By Proposition \ref{PsiConv_1}, $\bm{\Sigma}_{\Psi_n^{IID}({\bm{\theta}}_0)}- \bm{\Sigma}_{\Psi_n^{LHS}({\bm{\theta}}_0)}$ is asymptotically positive semi-definite, which implies the same property for $\bm{\Sigma}_{\bm{\hat{\theta}}_n^{IID}} - \bm{\Sigma}_{\bm{\hat{\theta}}_n^{LHS}}$.

Finally, if $\mathbb{E}(\|\psi_{\theta_0}(\bm{X})\|^3) < +\infty$ and $\bm{R}_{\psi_{\bm{\theta}_0}}$ is non-singular, the asymptotic normality follows from assertion 3 of Proposition \ref{PsiConv_1}. $\blacksquare$

These results give an asymptotic convergence for $\hat{\bm{\theta}}_n^{LHS}$ with, in the univariate case, a lower asymptotic variance of estimation than $\hat{\bm{\theta}}_n^{IID}$ (corresponding to $\bm{\Sigma}_{\bm{\hat{\theta}}_n^{IID}} - \bm{\Sigma}_{\bm{\hat{\theta}}_n^{LHS}}$ being asymptotically positive semi-definite in the multivariate case). Moreover, it gives a Central Limit Theorem for $Z$-estimators under LHS. Although strong regularity conditions on $\psi_{\bm{\theta}}$ are needed for these results to be valid, it remains very useful in many practical cases (e.g., for estimation by maximum likelihood). In the next section, we give an example of application. 

\section{Application: parameters estimation of Generalized Linear Models (GLM)} \label{Ex}

When performing statistical analysis on a computational code, it is common to approximate its outputs using a regression or classification model, also known as a metamodel. If the estimation of the modeling parameters can be expressed as a $Z$-estimator and the other conditions of use are satisfied, Theorem \ref{ZLHS_TCL} ensures that the estimation variance of these parameters is asymptotically lower under LHS than under IID sampling. It also provides a Central Limit Theorem under LHS.

Consider for instance the case of Generalized Linear Models (GLM), proposed in \cite{GLM}. They were formulated as a way of unifying various statistical models, including linear regression, logistic regression and Poisson regression. To estimate the parameters of a GLM,  one generally uses a Maximum Likelihood Estimator (MLE). It is therefore a special case of $Z$-estimation supposing that the likelihood can be differentiated. Thus, the results presented above can be applied to parameters estimation of a GLM.

\subsection{Definitions and main properties on GLM} \label{Ex1}

Before entering in more details, let us first define GLM more formally. For simplicity and without loss of generality, we focus here on the canonical case. Let $Z$ be a random variable on $\mathcal{Z} \subset \mathbb{R}$ and $\bm{X} = (X_1, \ldots, X_d)$ a vector of covariables on $\mathcal{X} \subset \mathbb{R}^d$, $d \in \mathbb{N}^{*}$. A GLM is characterized by:

\begin{enumerate}
    \item \textbf{A probability distribution:} $Z$ follows an exponential family distribution with density:
    \begin{equation}
        f(z,\alpha, \phi) = a(\alpha)b(z) \exp{(z \frac{\alpha}{\phi})}, \quad z \in \mathcal{Z}, \alpha \in \mathbb{R}, \phi > 0,
        \label{eq:density_GLM}
    \end{equation}
    where $\phi$ is the known dispersion parameter, $a(\alpha) = \exp(-v(\alpha)/\phi)$, $v: \mathbb{R} \to \mathbb{R}$ is twice continuously differentiable, and $b(z) = \exp(w(z, \phi))$ with $w: \mathcal{Z} \times \mathbb{R}^+ \to \mathbb{R}$ being also twice continuously differentiable in $z$.

    \item \textbf{A linear predictor:} For $\bm{\theta} = (\theta_1, \ldots, \theta_d)^T \in \Theta \subset \mathbb{R}^d$ , $\Theta$ being open and bounded:
    \begin{equation}
       \eta: \mathcal{X} \times \Theta \to \mathbb{R}, \quad \eta(\bm{x}, \bm{\theta}) = \bm{x}^T\bm{\theta} = \sum_{j=1}^d x_j \theta_j
       \label{eq:linear_GLM}.
    \end{equation}

    \item     \textbf{A link function:} Let $h: \mathcal{H} \subset \mathbb{R} \to \mathbb{R}$ be a monotone, differentiable function and $\bm{\alpha} = (\alpha^{(1)}, \ldots, \alpha^{(n)})^T \in \mathbb{R}^n$
    \begin{equation}
        h(\mu(\bm{x})) = \eta(\bm{x}, \bm{\theta}),
        \label{eq:link_function}
    \end{equation}
    where $\mu(\bm{x}) = \mathbb{E}[Z|\bm{X}=\bm{x}]$. Note that this hypothesis on $h$ implies the existence of the inverse function $h^{-1}: \mathbb{R} \rightarrow \mathcal{H}$ so that, for any $a \in \mathcal{H}$, $h^{-1} \circ h (a) = a$.
\end{enumerate}

Given $n$ independent realizations $\{(\bm{x}^{(i)}, z^{(i)})\}_{i=1}^n$, with $\bm{\textbf{X}} = (\bm{x}^{(1)},\ldots, \bm{x}^{(n)})^T$ and $\bm{\mathrm{Z}} = (z^{(1)},\ldots,z^{(n)})^T$, we also have $\bm{\alpha} = (\alpha^{(1)}, \ldots, \alpha^{(n)})^T \in \mathbb{R}^n$ such that:

\begin{equation}
    \alpha^{(i)} = h^{-1}(\bm{x}^{(i)T}\bm{\theta}), \quad i \in \llbracket 1,n \rrbracket.
\end{equation}

In this framework, the log-likelihood of each observation is:
\begin{equation}
    l(z^{(i)},\bm{x}^{(i)},\bm{\theta}, \phi) = \frac{1}{\phi}[z^{(i)}\alpha^{(i)}-v(\alpha^{(i)})] + w(z^{(i)}, \phi).
\end{equation}

The maximum likelihood estimator $\hat{\bm{\theta}}_n$ can be obtained by maximizing the log-likelihood function over the parameter space $\Theta$. Under regularity conditions, including the continuity and differentiability of the log-likelihood function over $\Theta$, the maximum likelihood estimator satisfies the first-order optimality conditions and therefore the following vectorial equation:

\begin{equation}
    \label{eq:ML_GLM}
    \sum_{i=1}^n \nabla_{\bm{\theta}} l(z^{(i)},\bm{x}^{(i)},\bm{\theta}, \phi) = \bm{0}.
\end{equation}

This defines $\hat{\bm{\theta}}_n$ as a $Z$-estimator with $\psi_{\bm{\theta}}(\bm{x}^{(i)}) = \nabla_{\bm{\theta}} l(z^{(i)},\bm{x}^{(i)},\bm{\theta}, \phi)$.

For the canonical case, the components of $\psi_{\bm{\theta}} = (\psi_{\theta_1}(\bm{x}), \ldots, \psi_{\theta_d}(\bm{x}))^T$ are:
\begin{equation}
    \psi_{\theta_j}(\bm{x}) = \frac{z - h^{-1}(\bm{x}\bm{\theta})}{\phi}x_j, \quad j \in \llbracket1, d\rrbracket, \bm{x} \in \mathcal{X}, \bm{z} \in \mathcal{Z}, \bm{\theta} \in \Theta.
\end{equation}
We can see that the estimation of the parameters of a GLM by maximum likelihood fits into the framework of $Z$-estimation. Thus, let us suppose that the observations of $\bm{X}$ are obtained by a LHS. We consider the $Z$-estimator defined by the equation \ref{eq:ML_GLM}, even though in this case the realizations are no longer \emph{i.i.d.} Let us discuss the convergence of this estimator under LHS.

\subsection{$Z$-estimation of GLM parameters under LHS} \label{Ex2}

Let $\mathrm{X}^{LHS} = (\bm{x}^{(i)}, \ldots \bm{x}^{(n)})^T$ be the realizations of $\bm{X}$ generated by a LHS. As before, for simplicity, we assume that we have $\mathcal{X} =[0,1]^d$ and $\bm{X} \sim U_{[0,1]^d}$. We also assume that $\Theta$ is the interior of a compact subset of $\mathbb{R}^d$. Moreover, we suppose that $h$ is defined and is twice continuously derivable on $\mathcal{H}$.  We also suppose that $h$ and its first derivative $\dot{h}: \mathcal{H} \to \mathbb{R}$ have no zero on $\mathcal{H}$. Since $h$ is monotone by construction, note that $h^{-1}$ is also defined and twice continuously derivable for any $\theta \in \Theta$ and $x \in \mathcal{X}$ thanks to the inverse function theorem (see, for instance, \cite{InversibleFunction} for more details).

Since we suppose that $\bm{X}$ and $\bm{\theta}$  are bounded, we have that $\Psi_n^{LHS}(\bm{\theta})=  \frac{1}{n}\sum_{i=1}^n \psi_{\bm{\theta}}(\bm{x}^{(i)})$ converges in probability to $\mathbb{E}(\psi_{\bm{\theta}}(\bm{X})) = \Psi(\bm{\theta})$. As we have seen, the other conditions concerning the convergence of ${\bm{\hat{\theta}}}_n$ to ${\bm{\theta}}_0$ are not specific to the sampling scheme. The conditions of application of Proposition \ref{convergence_2} are verified both in the case of an IID or a LHS design. We can thus conclude that ${\bm{\hat{\theta}}}_n$  converges in probability in ${\bm{\theta}}_0$.

Let us now verify that the conditions of application of Theorem \ref{ZLHS_TCL} are fulfilled. First, we see that $(\mathrm{X}^{LHS}, \bm{\theta}) \rightarrow \Psi_n^{LHS}(\bm{\theta})$ is continuous and twice continuously differentiable in $\bm{\theta}$. Plus, $\mathbb{E}(\psi_{\bm{\theta_0}}(\bm{X})) = \Psi(\bm{\theta_0}) = 0$ by construction.

We also have, for $j,k \in \llbracket 1,d \rrbracket$ , $\bm{x} =(x_1, \ldots, x_d)^T \in \mathcal{X}$ and $\bm{\theta} \in \Theta$:
\begin{equation}
 \label{eq:Psi_dot_GLM}
\frac{\partial \psi_{{\theta}_j}(\bm{x})}{\partial \theta_{k}} = \frac{-1}{\phi \dot{h}(h^{-1}(\bm{x} \bm{\theta}))} x_jx_{k}.
\end{equation}

Thus, we have that the matrix of partial derivatives $\bm{x}\rightarrow \dot{\psi}_{\bm{\theta}_0}(\bm{x})$ is such that $\mathbb{E}(\dot{\psi}_{\bm{\theta}_0}({X}))$ is defined and non-singular since $h$ has no zero on $\mathcal{H}$.  Since the values of $\bm{X}$ and $\bm{\theta}$ are bounded in norm, the function $(\bm{x}, \bm{\theta}) \rightarrow \psi_{\bm{\theta}}(\bm{x})$ is bounded and thus, for any $\bm{\theta} \in \Theta$,  $\mathbb E(||\psi_{\bm{\theta}}(X)||^3)<+\infty$ (and especially for $\bm{\theta} = \bm{\theta_0}$). 

Finally, we have that the elements of the tensor $\ddot{\psi}_{\bm{\theta}}(\bm{x})$ are, for  $j,k,l \in \rrbracket 1,d \llbracket$ and $\bm{x} =(x_1, \ldots, x_d)^T \in \mathcal{X}$, as follows:

\begin{equation}
 \label{eq:Psi_dot_dot_GLM}
\frac{\partial^2 \psi_{{\theta}_j}(\bm{x})}{\partial \theta_{k} \partial \theta_{l}} = \frac{\ddot{h}(h^{-1}(\bm{x} \bm{\theta}))}{\phi (\dot{h}(h^{-1}(\bm{x} \bm{\theta})))^3} x_jx_{k}x_{l}.
\end{equation}
Here $\ddot{h}: \mathcal{H}  \to \mathbb{R}$ is the second order derivative of $h$. 

Thus, $(\bm{x}, \bm{\theta})  \rightarrow ||\ddot{\psi}_{\bm{\theta}}(\bm{x})||$ can be bounded by an integrable function with finite second order moment in the neighborhood of $\bm{\theta_0}$ since we assume that the values of $\bm{\theta}$ and $\bm{X}$ are bounded.

All of these statements allow us to apply Theorem \ref{ZLHS_TCL}. We therefore have that the covariance matrix of estimation of $\hat{\bm{\theta}}_n^{LHS}$ is equal to $\bm{\Sigma}_{\bm{\hat{\theta}}_n^{LHS}} = [\mathbb{E}(\dot{\psi}_{\bm{\theta}_0}(\bm{X}))]^{-1}
\bm{\Sigma}_{\Psi_n^{LHS}({\bm{\theta}}_0)}[\mathbb{E}(\dot{\psi}_{\bm{\theta}_0}(\bm{X}))]^{-T} + \frac{1}{n}o(1)$ and that $\bm{\Sigma}_{\bm{\hat{\theta}}_n^{IID}} - \bm{\Sigma}_{\bm{\hat{\theta}}_n^{LHS}}$ is asymptotically positive semi-definite. Note that we have, with the previously introduced notations, $\bm{\Sigma}_{\Psi_n^{LHS}({\bm{\theta}}_0)} = \frac{1}{n}\bm{R}_{\psi_{\bm{\theta}_0}} + \frac{1}{n}o(1)$ asymptotically, with $\bm{R}_{\psi_{\bm{\theta}_0}} =\int_{[0,1]^d}{\psi_{{\bm{\theta}_0}_{rem}}}(\bm{x})\psi_{{\bm{\theta}_0}_{rem}}(\bm{x})^T\bm{dx}$.

Finally, since we have $\mathbb{E}(||\psi_{\theta_0}(\bm{X})||^3)<+\infty$, we have  that $\sqrt{n}(\hat{\bm{\theta}}_n^{LHS}- \bm{\theta}_0)$ is asymptotically normal with mean zero and a covariance matrix equal to $[\mathbb{E}(\dot{\psi}_{\bm{\theta}_0}(\bm{X}))]^{-1}\bm{R}_{\psi_{\bm{\theta}_0}}[\mathbb{E}(\dot{\psi_{\bm{\theta}_0}}(\bm{X}))]^{-1}$, assuming that $\bm{R}_{\psi_{\bm{\theta}_0}}$ is non-singular.

\subsection{Numerical example: a Poisson regression under LHS} \label{Ex3} 

To illustrate this result, let us consider a numerical example with a count random variable $Z$ and a vector of covariables $\bm{X} = (\bm{X}_1, \ldots, \bm{X}_{9})^T$, $\bm{X} \sim U_{[0,1]^{9}}$.  In an industrial context, $Z$ could represent for instance the number of operating problems evaluated by a simulation code of an industrial facility.  In this example, we define $Z$  by the following Poisson density function, for $z \in \mathbb{N}^{*}$:

\begin{equation}
    f(z, \lambda_0) = \exp(-\lambda_0)\frac{1}{z!} \exp{(z \log(\lambda_0))} ,
\end{equation}

with  $\log(\lambda_0) =  \bm{x}{\bm{\theta}_0}, \bm{x} \in [0,1]^9$ and $\bm{\theta}_0 = (\theta_{0,1}, \ldots \theta_{0,9})^T = (7,-\sqrt{2},1/2,-1/3,\sqrt{5},-7,\sqrt{2},-1/2,-\sqrt{5})^T$.

One can notice $Z$ fits in the framework of Equations \ref{eq:density_GLM}, \ref{eq:linear_GLM} and \ref{eq:link_function}.

Let us compare numerically the performances of the maximum likelihood estimation of $\bm{\theta}_0$ regarding the sampling method (IID or LHS) in this example. To do so, we first compare the estimation variance of each parameter $(\theta_{0,1}, \ldots, \theta_{0,9})^T$ with respect to the sampling scheme and size. We also verify that there is no significant difference concerning the square bias of estimation $\big[\mathbb{E}(\hat{\theta_j})-\theta_{0,j} \big]^2, j \in \llbracket 1,9 \rrbracket$. Additionally, we display the Mean Squared Error (MSE) for each parameter, defined as:

\begin{equation}
\text{MSE}(\hat{\theta_j}) = \mathbb{E}[(\hat{\theta_j} - \theta_{0,j})^2] = \text{Var}(\hat{\theta_j}) + [\mathbb{E}(\hat{\theta_j})-\theta_{0,j}]^2,
\end{equation}

where $\hat{\theta_j}$ is the estimator of $\theta_{0,j}$ and $j \in \llbracket 1,9 \rrbracket$. 

For each sampling method, the average values of these three metrics (variance, squared bias, and MSE) are computed over $L = 1000$ independent LHS and IID designs with sample sizes $n$ ranging from $40$ to $100$ (in increments of $10$).

\begin{figure}[ht!]
\centering
\includegraphics[width=0.55\textwidth]{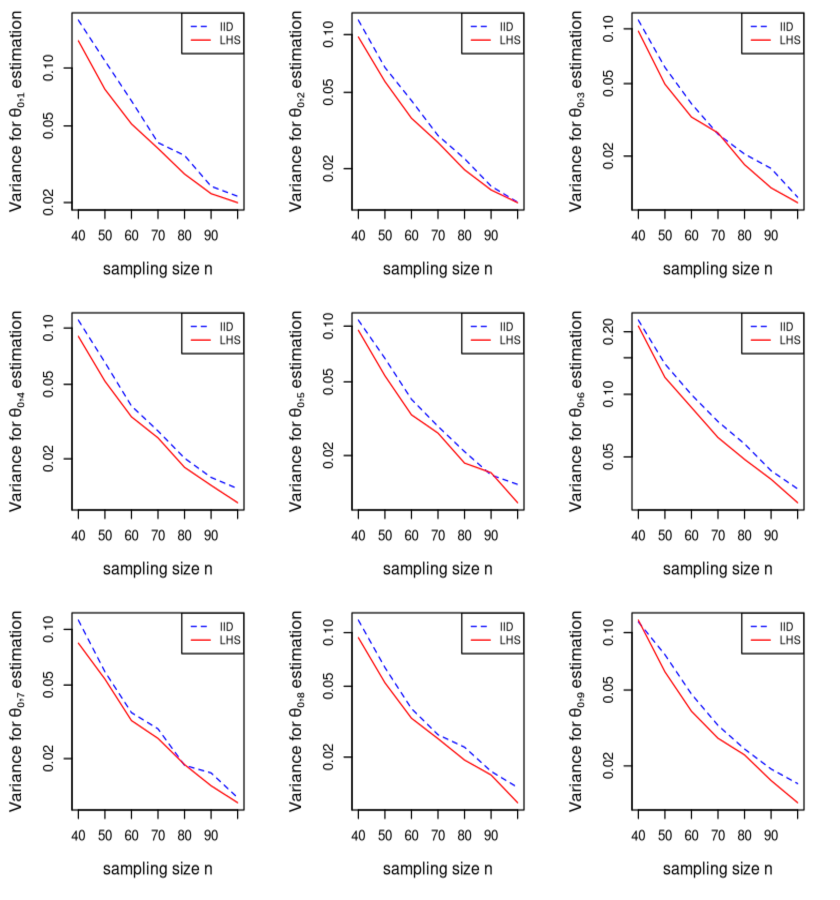}
\caption{Average estimation variances  of $(\theta_{0,1}, \ldots \theta_{0,9})^T$ according to the sampling size $n$ for IID and LHS designs (decimal logarithmic scale). \label{var_ex}}
\end{figure}

\begin{figure}[ht!]
\centering
\includegraphics[width=0.55\textwidth]{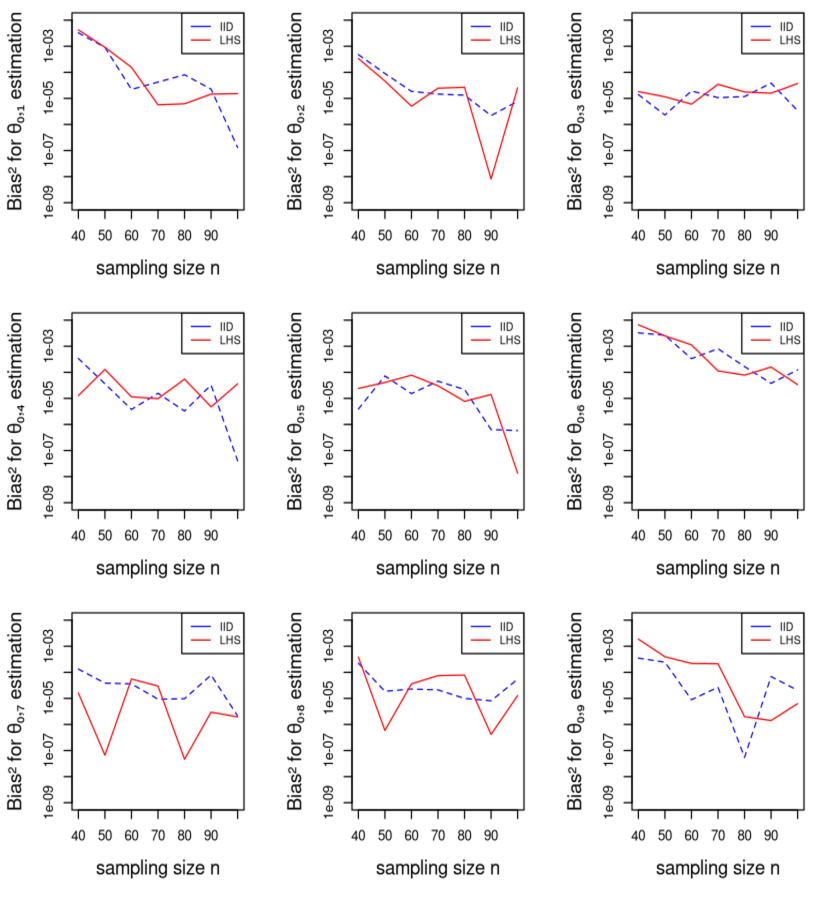}
\caption{Average estimation square bias  of $(\theta_{0,1}, \ldots \theta_{0,9})^T$ according to the sampling size $n$ for IID and LHS designs (decimal logarithmic scale). \label{bias_ex}}
\end{figure}

\begin{figure}[ht!]
\centering
\includegraphics[width=0.55\textwidth]{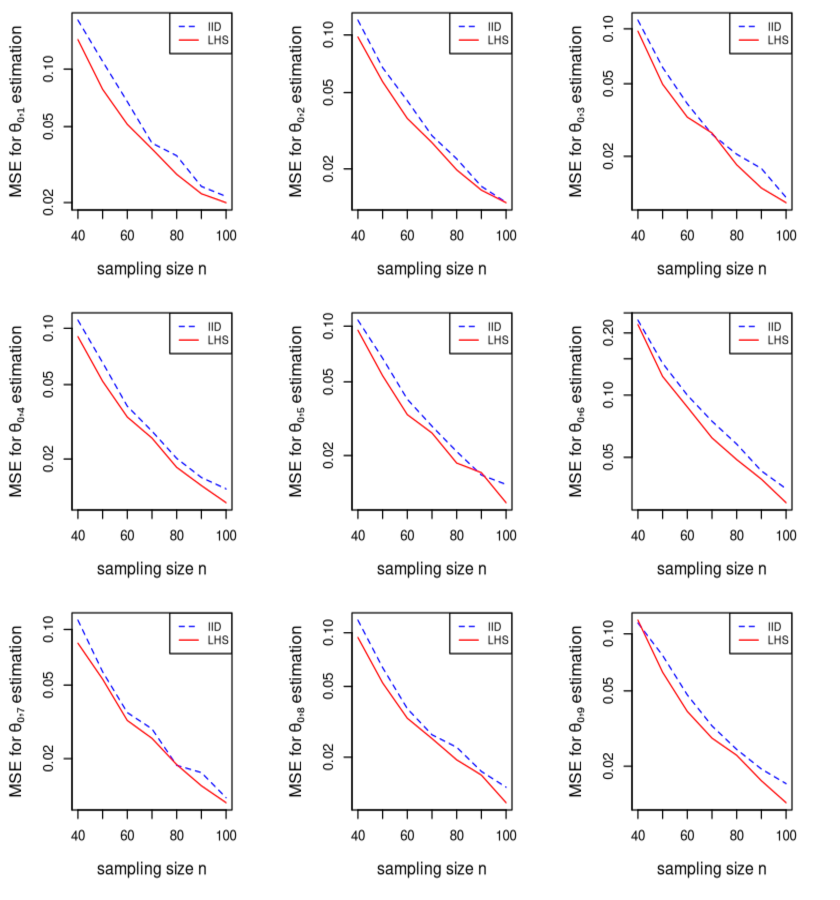}
\caption{Average estimation MSE  of $(\theta_{0,1}, \ldots \theta_{0,9})^T$ according to the sampling size $n$ for IID and LHS designs (decimal logarithmic scale). \label{MSE_ex}}
\end{figure}

Figure \ref{var_ex}, \ref{bias_ex} and \ref{MSE_ex} show respectively the evolution of the variance, the square bias of estimation and the MSE of the nine estimated parameters $(\theta_{0,1}, \ldots \theta_{0,9})^T$. As expected, we observe that for the nine estimated parameters, the average variance of estimation is overall lower for the classic LHS design compared to IID. No significant differences between LHS and IID designs are observed in terms of the square bias of estimation. The MSE is also significantly lower. As shown previously, classic LHS designs allow better estimation performances than IID ones, regardless of the theoretical value of the estimated parameters.

To further validate our theoretical results, we conduct an additional numerical experiment comparing the normalized empirical estimation variances ($n \times \text{Var}(\hat{\theta_j}), j \in \llbracket 1,9 \rrbracket$) with their theoretical asymptotic values. The covariance matrix formula, and hence the variance of each parameter (extracted from its diagonal), is derived from Theorem \ref{ZLHS_TCL}, using the application results given in \ref{Ex2} and Proposition \ref{var_mean}. The asymptotic variances evaluation is performed using a classic Monte Carlo method with a very large sample size ($n_{MonteCarlo} = 10^6$), as the analytical evaluation involves complex integrals and expectations that may not have closed-form solutions. The resulting theoretical asymptotic variances, rounded to four decimal places, are $(0.6195,0.3763,0.3388,0.3395,0.3963,1.4805,0.3575,0.3423,0.4328)$.

\begin{figure}[ht!]
\centering
\includegraphics[width=0.55\textwidth]{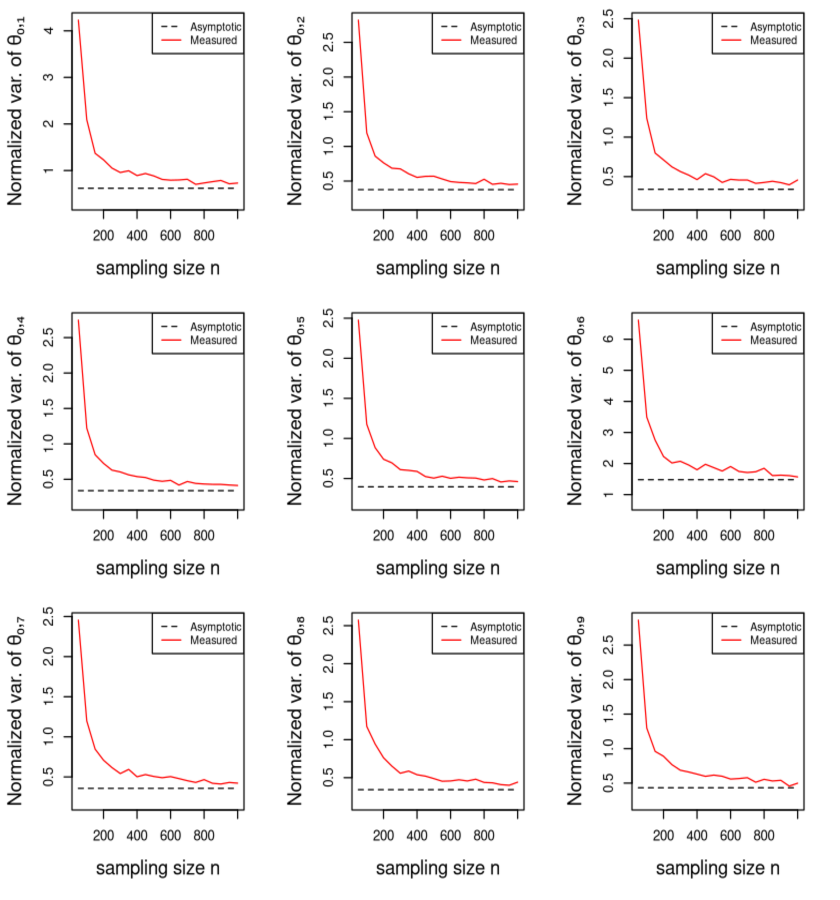}
\caption{Evolution of measured normalized variances of $\theta_{0,1}, \ldots \theta_{0,9}$ according to the sampling size $n$ for IID and LHS designs (decimal logarithmic scale). \label{var_conv_ex}}
\end{figure}

Figure \ref{var_conv_ex} shows the convergence of the normalized empirical variances with their theoretical counterparts as the sample size increases from $n = 50$ to $n=1000$ in increments of $50$. Each of the nine subplots demonstrates that the empirical normalized variances converge to the asymptotic values as $n$ grows, confirming the accuracy of our theoretical predictions. The average values are computed over $L =1000$ independent LHS designs, reinforcing the reliability of the observed trends.

\begin{figure}[ht!]
\centering
\includegraphics[width=0.55\textwidth]{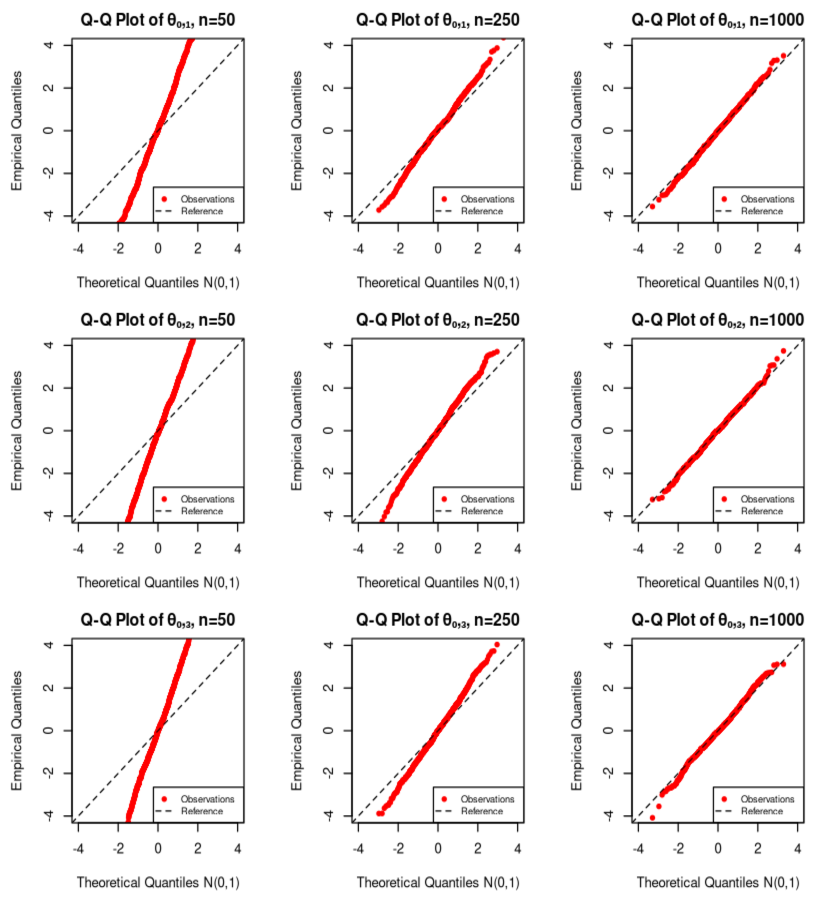}
\caption{Q-Q plots of the first three parameter estimates $\theta_{0,1}$, $\theta_{0,2}$, and $\theta_{0,3}$ (red dots) compared to the standard normal distribution for sample sizes of $50$, $250$, and $1000$. \label{qqplot_1}}
\end{figure}

\begin{figure}[ht!]
\centering
\includegraphics[width=0.55\textwidth]{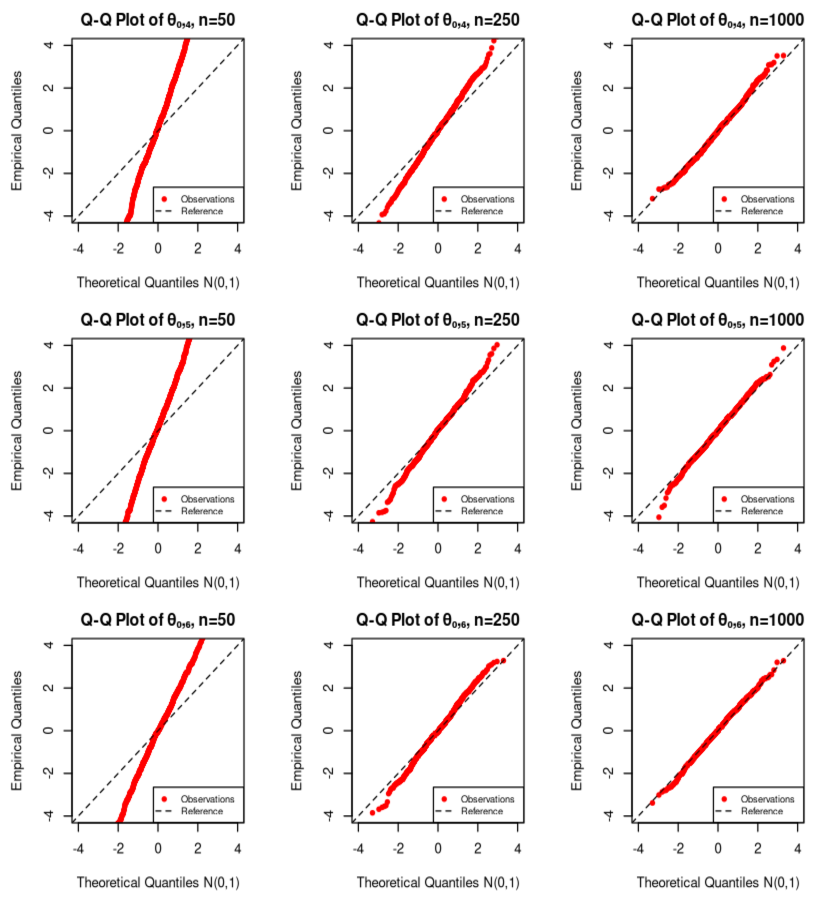}
\caption{Q-Q plots of the parameter estimates $\theta_{0,4}$, $\theta_{0,5}$, and $\theta_{0,6}$ compared to the standard normal distribution for sample sizes of $50$, $250$, and $1000$. \label{qqplot_2}}
\end{figure}

\begin{figure}[ht!]
\centering
\includegraphics[width=0.55\textwidth]{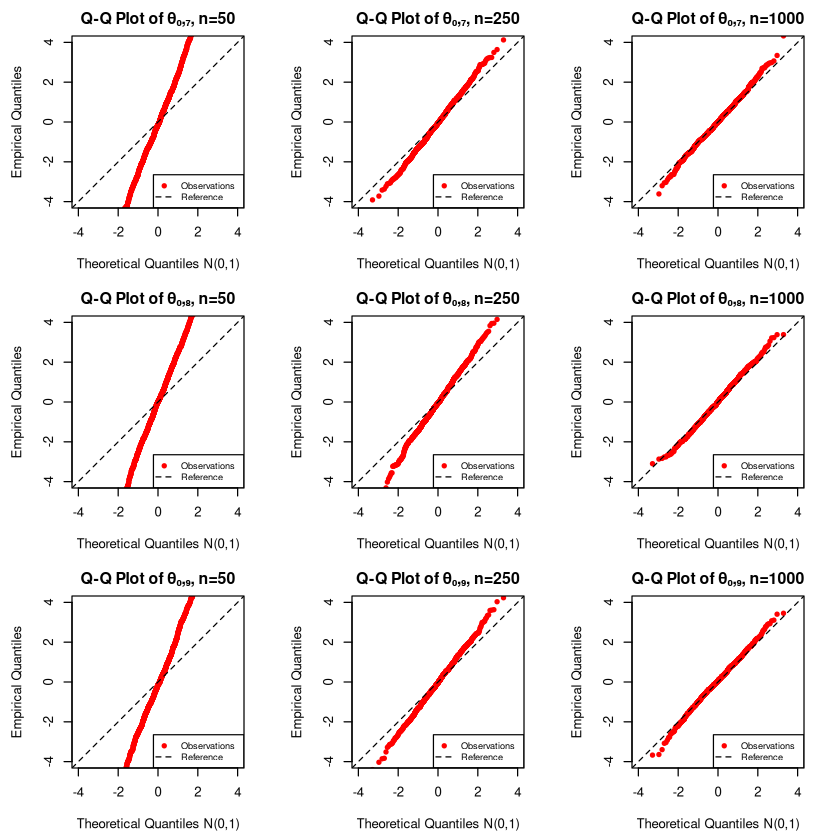}
\caption{Q-Q plots of the parameter estimates $\theta_{0,7}$, $\theta_{0,8}$, and $\theta_{0,9}$ compared to the standard normal distribution for sample sizes of $50$, $250$, and $1000$.  \label{qqplot_3}}
\end{figure}

Furthermore, to explore the asymptotic properties of the estimators, we conduct a third experiment examining the Q-Q plots (Quantile–Quantile plots) of the normalized parameter estimates compared to the standard normal distribution $\mathcal{N}(0,1)$. These plots compare the quantiles of the empirical distribution to the quantiles of a theoretical distribution, providing a visual assessment of how closely the data follows the specified distribution. Specifically, we generate these Q-Q plots for sample sizes of $n = 50$, $n = 250$, and $n = 1000$. Each plot is obtained with $L = 1000$ independent LHS designs. For clarity, we present the Q-Q plots for the first three parameters in Figure \ref{qqplot_1}, and similar plots for the remaining six parameters are provided in Figures \ref{qqplot_2} and \ref{qqplot_3}. The Q-Q plots reveal that the empirical distribution exhibits heavy tails for $n=50$, which improve with $n=250$ and align even more closely with the theoretical distribution for $n=1000$. This comparison allows us to visually assess the convergence of the estimators to normality, as predicted by our theoretical results.

In summary, our numerical experiments confirm that Latin Hypercube Sampling (LHS) consistently outperforms Independent and Identically Distributed (IID) sampling in terms of estimation variance and Mean Squared Error (MSE). The normalized empirical variances converge to their theoretical values as the sample size increases, validating our theoretical work. Additionally, the Q-Q plots visually confirm the asymptotic normality of the estimators.

\hypertarget{var_ex}{}

\hypertarget{bias_ex}{}

\hypertarget{MSE_ex}{}

\hypertarget{var_conv_ex}{}

\hypertarget{qqplot_1}{}

\hypertarget{qqplot_2}{}

\clearpage

\hypertarget{qqplot_3}{}


\section{Conclusion and prospects} \label{Cl}

In conclusion, Latin Hypercube Sampling (LHS) has demonstrated its robustness as a powerful method for conducting computer experiments, particularly in the analysis of complex black-box functions. This paper has advanced the understanding of the asymptotic convergence of estimators using this sampling method. Specifically, we have extended the convergence results under LHS previously established in \cite{Stein}, \cite{Owen}, and \cite{Loh} for the empirical mean to the broader class of $Z$-estimators. A key contribution of this work is the introduction of a Central Limit Theorem (CLT) for $Z$-estimators under this sampling method with a reduced asymptotic variance compared to traditional independent and identically distributed (i.i.d.) random sampling. Furthermore, we have illustrated the practical relevance of these theoretical findings through an application to parameter estimation in Generalized Linear Models (GLMs) under LHS. However, it is worth noting that certain restrictive regularity conditions were necessary to establish these convergence results.

A promising direction for future research is to relax some of the regularity assumptions, particularly the requirement for the second derivative of the $Z$-function. Alternative formulations of the CLT for $Z$-estimators  that do not depend on the existence of a second derivative have been proposed, as discussed in \cite{VanderVaart}. Exploring these approaches for the convergence of $Z$-estimators under LHS could lead to more generalized results that extend beyond those presented in this work. Furthermore, generalizing the given convergence results under LHS to the class of $M$-estimators presents another valuable perspective for future research.

Another valuable extension of this research would be to adapt our theoretical framework to models with complex dependency structures, such as Gaussian Process (GP) regression \cite{Gaussian_Process,Kang_GP,Marrel_GP}, a key application of LHS. While GP regression parameter estimation does not align with the $Z$-estimator framework due to its non-separable log-likelihood function, developing methodologies to address these dependencies could substantially expand the scope and applicability of our findings.

Ultimately, this work underscores the significant value of LHS in industrial applications, particularly for analyzing simulation codes that are computationally intensive and involve numerous input parameters. The versatility of LHS enables the efficient implementation of various statistical techniques—including variable selection, sensitivity analysis, and metamodeling—within a single numerical design of experiments. For further exploration of practical industrial applications of LHS, we refer readers to \cite{Hakimi, Aleksander, BOURCET, FERRARI, Marrel2} as a few examples among many available in the literature.

\section*{Acknowledgments}

Support from my supervisors (Claude Brayer, Fabrice Gamboa, Benoît Habert and Amandine Marrel) during my PhD thesis, of which this work is a continuation, is gratefully acknowledged. I would also like to thank Anouar Meynaoui for his help in writing this paper.



\bibliography{mybib}

\begin{thebibliography}{31}
\providecommand{\natexlab}[1]{#1}
\providecommand{\url}[1]{\texttt{#1}}
\providecommand{\urlprefix}{}

\bibitem[{Mckay et~al.(1979)Mckay, M. and Beckman, Richard and Conover, William}]{LHS}
Mckay M, Beckman R, Conover W.
\newblock A Comparison of Three Methods for Selecting Vales of Input Variables in the Analysis of Output From a Computer Code 21.
\newblock Technometrics 1979;p. 239--245.

\bibitem[{Helton and Davis(2003)J.C. Helton and F.J. Davis}]{Helton}
Helton JC, Davis FJ.
\newblock Latin hypercube sampling and the propagation of uncertainty in analyses of complex systems.
\newblock Reliability Engineering \& System Safety 81 2003;p. 23--69.

\bibitem[{Helton et~al.(2005)J.C. Helton and F.J. Davis and J.D. Johnson}]{Helton2}
Helton JC, Davis FJ, Johnson JD.
\newblock A comparison of uncertainty and sensitivity analysis results obtained with random and Latin hypercube sampling.
\newblock Reliability Engineering \& System Safety 89 2005;p. 305--330.

\bibitem[{Viana(2016)Viana, Felipe AC}]{viana2016}
Viana FA.
\newblock A tutorial on Latin hypercube design of experiments.
\newblock Quality and reliability engineering international 2016;32(5):1975--1985.

\bibitem[{Deutsch and Deutsch(2012)Jared L. Deutsch and Clayton V. Deutsch}]{DEUTSCH}
Deutsch JL, Deutsch CV.
\newblock Latin hypercube sampling with multidimensional uniformity.
\newblock Journal of Statistical Planning and Inference 2012;142(3):763--772.

\bibitem[{Shields and Zhang(2016)Michael D. Shields and Jiaxin Zhang}]{SHIELDS}
Shields MD, Zhang J.
\newblock The generalization of Latin hypercube sampling.
\newblock Reliability Engineering \& System Safety 2016;148:96--108.

\bibitem[{Packham and Schmidt(2008)Packham, Natalie and Schmidt, Wolfgang M}]{packham2008latin}
Packham N, Schmidt WM.
\newblock Latin hypercube sampling with dependence and applications in finance.
\newblock Available at SSRN 1269633 2008;.

\bibitem[{Mondal and Mandal(2020)Anirban Mondal and Abhijit Mandal}]{MONDAL}
Mondal A, Mandal A.
\newblock Stratified random sampling for dependent inputs in Monte Carlo simulations from computer experiments.
\newblock Journal of Statistical Planning and Inference 2020;205:269--282.

\bibitem[{Stein(1987)Stein, Michael}]{Stein}
Stein M.
\newblock Large {Sample} {Properties} of {Simulations} {Using} {Latin} {Hypercube} {Sampling}.
\newblock Technometrics 2 1987;p. 143--151.

\bibitem[{Owen(1992)Owen, Art B.}]{Owen}
Owen AB.
\newblock A {Central} {Limit} {Theorem} for {Latin} {Hypercube} {Sampling}.
\newblock Journal of the Royal Statistical Society 54: Series B (Methodological) 1992;p. 541--551.

\bibitem[{Loh(1996)Loh, Wei Liem}]{Loh}
Loh WL.
\newblock On latin hypercube sampling.
\newblock Annals of Statistics 24 1996;p. 2058–2080.

\bibitem[{Aistleitner et~al.(2013)Christoph Aistleitner and Markus Hofer and Robert F. Tichy}]{Aistleitner2013}
Aistleitner C, Hofer M, Tichy RF.
\newblock A central limit theorem for Latin hypercube sampling with dependence and application to exotic basket option pricing.
\newblock International Journal of Theoretical and Applied Finance 2013;15(07).

\bibitem[{Hakimi(2023)Hakimi, Faouzi}]{Hakimi}
Hakimi F.
\newblock High-dimensional sensitivity analysis methods for computationally expensive simulators modeling a severe nuclear accident.
\newblock PhD thesis, Paul Sabatier University, Toulouse, France; 2023.

\bibitem[{Devroye(1986)Devroye, Luc}]{Tinverse}
Devroye L.
\newblock 2: General principles in random variate generation.
\newblock In: Non-Uniform Random Variate Generation(originally published with Springer-Verlag; 1986. p. 27--39.

\bibitem[{Rousseeuw et~al.(1986)Rousseeuw, P.J. and Hampel, F.R. and Ronchetti, E.M. and Stahel, W.A.}]{Robust}
Rousseeuw PJ, Hampel FR, Ronchetti EM, Stahel WA.
\newblock Robust Statistics: The Approach Based on Influence Functions.
\newblock Wiley; 1986.

\bibitem[{{Van der Vaart}(1998){Van der Vaart}, A. W.}]{VanderVaart}
{Van der Vaart} AW.
\newblock Asymptotic statistics.
\newblock Cambridge University Press; 1998.

\bibitem[{Bouzebda and Chaouch(2022)Salim Bouzebda and Mohamed Chaouch}]{Bouzebda2022}
Bouzebda S, Chaouch M.
\newblock Uniform limit theorems for a class of conditional Z-estimators when covariates are functions.
\newblock Journal of Multivariate Analysis 2022;189:104872.

\bibitem[{Zhan(2002)Yihui Zhan}]{Zhan2002}
Zhan Y.
\newblock Central limit theorems for functional Z-estimators.
\newblock Statistica Sinica 2002;12:609--634.

\bibitem[{Hall(2003)Hall, Alastair R}]{GMM}
Hall AR.
\newblock Generalized method of moments.
\newblock A companion to theoretical econometrics 2003;p. 230--255.

\bibitem[{Dacunha-Castelle and Duflo(1986)Dacunha-Castelle, D. and Duflo, M.}]{Duflo}
Dacunha-Castelle D, Duflo M.
\newblock Probability and Statistics Volume II.
\newblock Springer-Verlag; 1986.

\bibitem[{Fisher(????)Fisher, R.A.}]{Fisher_likelihood}
Fisher RA.
\newblock On the mathematical foundations of theoretical statistics.
\newblock Philosophical Transactions of the Royal Society of London, Series A 222;p. 309–368.

\bibitem[{Fischer(2011)Fischer, Hans}]{ClassicTCL}
Fischer H.
\newblock A History of the Central Limit Theorem: From Classical to Modern Probability Theory.
\newblock New York: Springer; 2011.

\bibitem[{Nelder and Wedderburn(????)Nelder, John and Wedderburn, Robert}]{GLM}
Nelder J, Wedderburn R.
\newblock Generalized Linear Models.
\newblock Journal of the Royal Statistical Society Series A (General) Blackwell Publishing 135;p. 370–384.

\bibitem[{Allendoerfer(1974)Allendoerfer, Carl B.}]{InversibleFunction}
Allendoerfer CB.
\newblock Calculus of Several Variables and Differentiable Manifolds.
\newblock New York: Macmillan; 1974.

\bibitem[{Williams and Rasmussen(1995)Williams, Christopher and Rasmussen, Carl}]{Gaussian_Process}
Williams C, Rasmussen C.
\newblock Gaussian processes for regression.
\newblock Advances in neural information processing systems 1995;8.

\bibitem[{Kang et~al.(2015)Kang, Fei and Han, Shaoxuan and Salgado, Rodrigo and Li, Junjie}]{Kang_GP}
Kang F, Han S, Salgado R, Li J.
\newblock System probabilistic stability analysis of soil slopes using Gaussian process regression with Latin hypercube sampling.
\newblock Computers and geotechnics 2015;63:13--25.

\bibitem[{Marrel and Iooss(2024)Marrel, Amandine and Iooss, Bertrand}]{Marrel_GP}
Marrel A, Iooss B.
\newblock Probabilistic surrogate modeling by Gaussian process: A review on recent insights in estimation and validation.
\newblock Reliability Engineering \& System Safety 2024;247:110094.

\bibitem[{Karolczuk and Kurek(2022)Aleksander Karolczuk and Marta Kurek}]{Aleksander}
Karolczuk A, Kurek M.
\newblock Fatigue life uncertainty prediction using the Monte Carlo and Latin hypercube sampling techniques under uniaxial and multiaxial cyclic loading.
\newblock International Journal of Fatigue 2022;160:106867.

\bibitem[{Bourcet et~al.(2023)J. Bourcet and A. Kubilay and D. Derome and J. Carmeliet}]{BOURCET}
Bourcet J, Kubilay A, Derome D, Carmeliet J.
\newblock Representative meteorological data for long-term wind-driven rain obtained from Latin Hypercube Sampling – Application to impact analysis of climate change.
\newblock Building and Environment 2023;228:109875.

\bibitem[{Ferrari et~al.(2019)Rosalba Ferrari and Diego Froio and Egidio Rizzi and Carmelo Gentile and Eleni N. Chatzi}]{FERRARI}
Ferrari R, Froio D, Rizzi E, Gentile C, Chatzi EN.
\newblock Model updating of a historic concrete bridge by sensitivity- and global optimization-based Latin Hypercube Sampling.
\newblock Engineering Structures 2019;179:139--160.

\bibitem[{Marrel et~al.(2022)Marrel, Amandine and Iooss, Bertrand and Chabridon, Vincent}]{Marrel2}
Marrel A, Iooss B, Chabridon V.
\newblock The icscream methodology: Identification of penalizing configurations in computer experiments using screening and metamodel—applications in thermal hydraulics.
\newblock Nuclear Science and Engineering 2022;196(3):301--321.

\end{thebibliography}

\end{document}